\documentclass[journal]{IEEEtran}

\makeatletter
\def\ps@headings{%
\def\@oddhead{\mbox{}\scriptsize\rightmark \hfil \thepage}%
\def\@ adversarynhead{\scriptsize\thepage \hfil \leftmark\mbox{}}%
\def\@oddfoot{}%
\def\@ adversarynfoot{}}
\makeatother
\ifCLASSINFOpdf
\else

\fi
\usepackage[linesnumbered,ruled,vlined]{algorithm2e}
\usepackage{algorithm2e}
\usepackage{epsfig}
\usepackage{array}
\newcolumntype{L}[1]{>{\raggedright\let\newline\\\arraybackslash\hspace{0pt}}m{#1}}
\newcolumntype{C}[1]{>{\centering\let\newline\\\arraybackslash\hspace{0pt}}m{#1}}
\newcolumntype{R}[1]{>{\raggedleft\let\newline\\\arraybackslash\hspace{0pt}}m{#1}}
\usepackage{amsmath}
\usepackage{amsthm} % enhanced theorem-like environments
\usepackage{mathrsfs}
\newcommand{\bc}{\begin{center}}
\newcommand{\ec}{\end{center}}
\newcommand{\be}{\begin{equation}}
\newcommand{\ee}{\end{equation}}
\usepackage{flexisym}
\usepackage{algorithmicx}
\usepackage{wasysym}%
\usepackage{soul}

\usepackage{amsmath}

\newcommand{\bnu}{\begin{enumerate}}
\newcommand{\enu}{\end{enumerate}}
\usepackage{cite}
\usepackage{url}
\usepackage{breqn}
\usepackage{lipsum}
\usepackage{mathtools}
\usepackage{amsmath}
\usepackage{caption}
\usepackage{subcaption}
\usepackage{comment}
\usepackage{soul}
\usepackage{blindtext, graphicx}
\usepackage{cuted}
\usepackage{color}

\usepackage{amsfonts}
\usepackage{lipsum}
\usepackage{cuted}
\usepackage{float}

\begin{document}
\pagestyle{empty}
%\newgeometry{top=1in,bottom=1.05in,right=0.75in,left=0.75in}
\title{QSAFE-V: Quantum-Enhanced Lightweight Authentication Protocol Design for Vehicular Tactile Wireless Networks}
%\author{ Authors

\author{Shakil Ahmed,~\IEEEmembership{Member,~IEEE}, Amika Tabassum, Ibrahim Almazyad, and Ashfaq Khokhar,~\IEEEmembership{Fellow,~IEEE}
\vspace*{-1.05 cm}
\thanks{ 
}
\thanks{ Shakil Ahmed, Amika Tabassum, and Ashfaq Khokhar are with the Department of Electrical and Computer Engineering, Iowa State University, Ames, Iowa, USA. (email: \{shakil, amika, ashfaq\}@iastate.edu) \\
I. Almazyad is with the Department of Computer Engineering, Al-Qassim University, KSA. (email: i.almazyad@qu.edu.sa). 
%Nasser Albalawi is with the Department of Electrical and Computer Engineering, Ozyegin University, 34794 ˙Istanbul, Türkiye (e-mail: Nasser.Albalawi@nbu.edu.sa).
}
}
\markboth{IEEE Transactions on Wireless Communications}
{}
 
\maketitle
 
 \thispagestyle{empty}
 
\IEEEpeerreviewmaketitle

\begin{abstract}
With the rapid advancement of 6G technology, the Tactile Internet is emerging as a novel paradigm of interaction, particularly in intelligent transportation systems, where stringent demands for ultra-low latency and high reliability are prevalent. During the transmission and coordination of autonomous vehicles, malicious adversaries may attempt to compromise control commands or swarm behavior, posing severe threats to road safety and vehicular intelligence. Many existing authentication schemes claim to provide security against conventional attacks. However, recent developments in quantum computing have revealed critical vulnerabilities in these schemes, particularly under quantum-enabled adversarial models. In this context, the design of a quantum-secured, lightweight authentication scheme that is adaptable to vehicular mobility becomes essential. This paper proposes QSAFE-V, a quantum-secured authentication framework for edge-enabled vehicles that surpasses traditional security models. We conduct formal security proofs based on quantum key distribution and quantum adversary models, and also perform context-driven reauthentication analysis based on vehicular behavior. The output of quantum resilience evaluations indicates that QSAFE-V provides robust protection against quantum and contextual attacks. Furthermore, detailed performance analysis reveals that QSAFE-V achieves comparable communication and computation costs to classical schemes, while offering significantly stronger security guarantees under wireless Tactile Internet conditions.
\end{abstract}


\section{Introduction}
\IEEEPARstart{W}{ith} the proliferation of ultra-reliable low-latency communication (URLLC) and advanced edge computing infrastructures, the vision of the Tactile Internet is rapidly evolving toward realization, which enables real-time remote control over both physical and virtual environments by transmitting haptic and motion data over the network, empowering applications such as telesurgery, collaborative robotics, and immersive touch-based interfaces~\cite{ngmn2015ti}.
However, as Tactile Internet becomes increasingly embedded in mission-critical healthcare applications, the need for secure and provable authentication becomes urgent. Lightweight cryptographic solutions, while efficient, are often vulnerable to sophisticated attacks—especially under quantum computing capabilities~\cite{chen2016pqcrypto}. Current schemes fail to provide security guarantees against adversaries equipped with quantum resources capable of breaking conventional hardness assumptions \cite{karobi2025ecoedgetwin}.

To address these challenges, we propose a Quantum Secured Authentication Framework for edge-enabled vehicles (QSAFE-V), which combines post-quantum cryptographic primitives with Quantum Key Distribution (QKD), providing unconditional security grounded in quantum mechanics~\cite{bennett1984qkd, pirandola2020advances}. 
%Additionally, the use of Physically Unclonable Functions (PUFs) offers strong protection against physical tampering and cloning, making the scheme ideal for secure medical edge devices~\cite{maes2010puf}.
This paper provides a comprehensive evaluation of QSAFE-V through both formal security analysis using an extended Real-Or-Random (ROR) model in vehicular networks. The proposed scheme demonstrates robustness against passive, active, and implicit attacks, ensuring confidentiality, integrity, anonymity, and resistance to desynchronization, even under post-quantum threat model for vehicular tactile wireless networks.
Moreover, QSAFE-V achieves its security guarantees with low communication and computation overhead, making it well-suited for deployment in URLLC-sensitive, resource-constrained environments like the Tactile Internet~\cite{itu2014tistandards, ngmn2015ti}.

%\subsection{Related Work}
Tactile Internet represents a transformative approach to real-time interaction over networks with applications ranging from remote surgery to autonomous vehicles \cite{ahmed2025quantum}. Authors in~\cite{fettweis2014tactile} introduced the Tactile Internet concept, highlighting URLLC as its cornerstone. The expansion into vehicular domains has been explored in~\cite{aijaz2016towards}, who presented the integration of haptic communication with 5G for vehicular control. The merging of Tactile Internet and Vehicular Edge Computing (VEC) has since been extended in~\cite{singh2021latencyaware, alsabah2021lightweight} to support intelligent traffic and remote vehicle control with latency constraints below 1 ms.
As classical cryptographic methods face obsolescence under quantum attack models, QKD and post-quantum cryptography (PQC) gained significant attraction. The authors in~\cite{kim2021vision} envisioned the role of QKD in future network architectures. The authors in~\cite{sun2022quantum} introduced a QKD-based authentication for healthcare tactile systems. In vehicular contexts, the authors in~\cite{misra2022leveraging} proposed QKD-resilient protocols for dynamic networks. Additionally, robust PQC schemes such as lattice-based and hash-based cryptography have been recommended for integration into vehicular environments~\cite{chen2016report,najafi2022towards}. Moreover, a Quantum Identity Token (QIT) is crucial for securing vehicular tactile networks by providing an unforgeable and tamper-evident mechanism to authenticate vehicles and control signals in real-time, which is essential for safety-critical applications (see Section~\ref{Sec_prelem}) for more details.

Quantum Machine Learning (QML) enhances processing speed and energy efficiency, vital for  Tactile Internet. The authors in~\cite{schuld2014quest} outlined quantum neural models with lower complexity. The authors in~\cite{biamonte2017quantum} showed how quantum states can accelerate learning processes. The hybrid model combining classical and quantum inference at edge nodes was introduced in~\cite{he2022hybrid}, demonstrating improvements in decision latency.
Several authentication frameworks for the Tactile Internet have been proposed, such as IAR-AKA~\cite{alsabah2021lightweight}, which uses hash and ECC-based lightweight primitives. However, these schemes remain vulnerable to quantum attacks. 
QKD has shown promise for mutual authentication with low overhead~\cite{alkhateeb2022efficient}. The authors in~\cite{pirandola2020advances} further elaborated on the need for quantum-secure network protocols for tactile use cases.

Classical authentication mechanisms, even when optimized for lightweight performance, remain inadequate in the face of emerging quantum threats and dynamic vehicular edge conditions.
Despite these advancements, several challenges remain: (1) most Tactile Internet authentication schemes do not account for quantum threat models, (2) hybrid quantum-classical frameworks are still experimental, and (3) quantum edge deployment in vehicular networks remains a bottleneck. To address these concerns, this paper introduces QSAFE-V that operates under the Tactile Internet paradigm.
 To the best of our knowledge, this is the first authentication framework that combines the principles of QIT, Tactile Internet design constraints, and VEN in a unified security protocol.  The key contributions of this work are as follows:

\begin{enumerate}

\item \textit{Quantum-Resilient Authentication Design:} We propose a novel authentication framework that integrates QKD with hash-based verification to achieve mutual authentication and forward secrecy, ensuring resilience against both classical and quantum-based attacks, such as replay, impersonation, and quantum brute-force.

\item \textit{Tactile Vehicular Edge Integration:} QSAFE-V is designed specifically for VEC environments, addressing the stringent latency and reliability requirements of Tactile Internet of Vehicles. 

\item \textit{Lightweight and Scalable Protocol:} Unlike conventional quantum protocols that assume heavy computation, QSAFE-V adopts hybrid classical-quantum communication primitives and minimizes quantum overhead. This enables deployment in resource-constrained edge nodes and supports horizontal scaling across vehicular networks.

\item \textit{Security and Performance Evaluation:} We perform extensive analytical evaluations to verify the proposed protocol's resistance to common attack vectors. Our results also demonstrate improvements in handshake latency and entropy strength compared to existing classical schemes such as IAR-AKA \cite{yang2025iar}.

\item \textit{Comprehensive Use Case Targeting:} The framework is applicable across a wide range of vehicular
Tactile wireless networks use cases, including remote tactile feedback systems, haptic-enabled vehicle control, and authentication with QIT-layer enhancements.

\end{enumerate}

%\section{Related Work}
%\input{sections/related_work}

\section{Preliminaries} \label{Sec_prelem}
QIT is a security primitive rooted in quantum information theory that enables device authentication through non-clonable quantum states. QITs exploit the quantum no-cloning theorem, which ensures that arbitrary quantum states cannot be perfectly duplicated. These tokens are generated using entangled photon pairs or quantum states with device-specific parameters and may be used to perform secure identification of Autonomous
Vehicles (AVs) within Tactile Internet environments.
Each QIT is represented as a tuple of a quantum challenge and its expected quantum response, governed by quantum measurement principles. Formally, the authentication process involves sending a quantum challenge $Q$ and receiving a quantum response $R = \text{Measure}(Q)$ under the device's internal quantum state.
QITs possess the following characteristics:
\begin{itemize}
    \item \textit{Quantum Unclonability:} The quantum state associated with each QIT cannot be cloned due to the no-cloning theorem. This makes it inherently secure against replication and impersonation attacks.
    
    \item \textit{Measurement-Driven Uniqueness:} The response to a quantum challenge depends on the internal quantum state of the device. Even slight deviations in device parameters will produce statistically distinct outcomes.
    
    \item \textit{Non-Predictability:} The response of a device to a quantum challenge is governed by probabilistic measurement outcomes, making it computationally infeasible to predict responses using classical or quantum computation in polynomial time.
\end{itemize}

\subsection{Post-Quantum Cryptography}
\subsubsection{Lattice-Based Cryptography}
In the context of post-quantum cryptography, lattice-based schemes are considered among the most promising and efficient candidates to replace classical elliptic curve cryptography. A lattice is defined as a discrete, periodic arrangement of points in $n$-dimensional space, generated by linear combinations of basis vectors with integer coefficients. The fundamental problems on which lattice cryptography relies are the Shortest Vector Problem (SVP) and the Learning With Errors (LWE) problem.
Lattice structures exhibit the following key properties:
\begin{itemize}
    \item \textit{Quantum Hardness:} The SVP and LWE problems are believed to be resistant to both classical and quantum attacks. No efficient quantum algorithm has been found to solve them in polynomial time.
    \item \textit{Noise-Based Security:} In LWE, small random errors are introduced during computation, making it computationally infeasible for an adversary to recover the original values even with quantum capabilities.
    \item \textit{Additive Homomorphism:} Certain lattice-based schemes support homomorphic operations, which can be leveraged for privacy-preserving authentication and lightweight computation at vehicular edge nodes.
\end{itemize}

\subsubsection{Hard Problems Under Lattices}
Lattice-based cryptography derives its security from two widely studied problems:
Given a secret vector $\mathbf{s} \in \mathbb{Z}_q^n$ and a public matrix $\mathbf{A} \in \mathbb{Z}_q^{m \times n}$, the LWE problem asks the adversary to distinguish between the noisy inner product $\mathbf{A} \cdot \mathbf{s} + \mathbf{e}$ (where $\mathbf{e}$ is a small error vector) and a uniform random vector in $\mathbb{Z}_q^m$. Solving this in polynomial time is considered infeasible even for quantum computers.
Given a matrix $\mathbf{A} \in \mathbb{Z}_q^{m \times n}$, the Short Integer Solution (SIS) problem asks to find a non-zero integer vector $\mathbf{x} \in \mathbb{Z}^n$ such that $\mathbf{A} \cdot \mathbf{x} = 0 \mod q$ and $\|\mathbf{x}\|$ is small. This forms the basis for secure digital signatures and commitments in quantum-resistant schemes.
In this work, the QSAFE-V protocol utilizes lattice-based cryptographic functions for key exchange, identity token encryption, and session establishment. These primitives ensure resistance against Shor's and Grover’s quantum attacks, making the scheme suitable for deployment in post-quantum wireless vehicular networks.

\section{System Model}
\subsection{Network Model}
The network model of QSAFE-V is shown in Fig.~\ref{fig:SM}. In the proposed model, there are primarily four entities: AVs, Vehicular Edge Nodes (VENs), Gateways (GW), and Trusted Authorities (TAs). Below, these entities in Tactile Internet environment are elaborately discussed:
\begin{comment}
    \begin{figure}[H]
\centering
\includegraphics[width=0.95\linewidth]{SystemModel.png}
\caption{QSAFE-V system model}
\label{fig:SM}
\end{figure}
\end{comment}
\begin{figure}[H]
\centering
\includegraphics[width=0.80\linewidth]{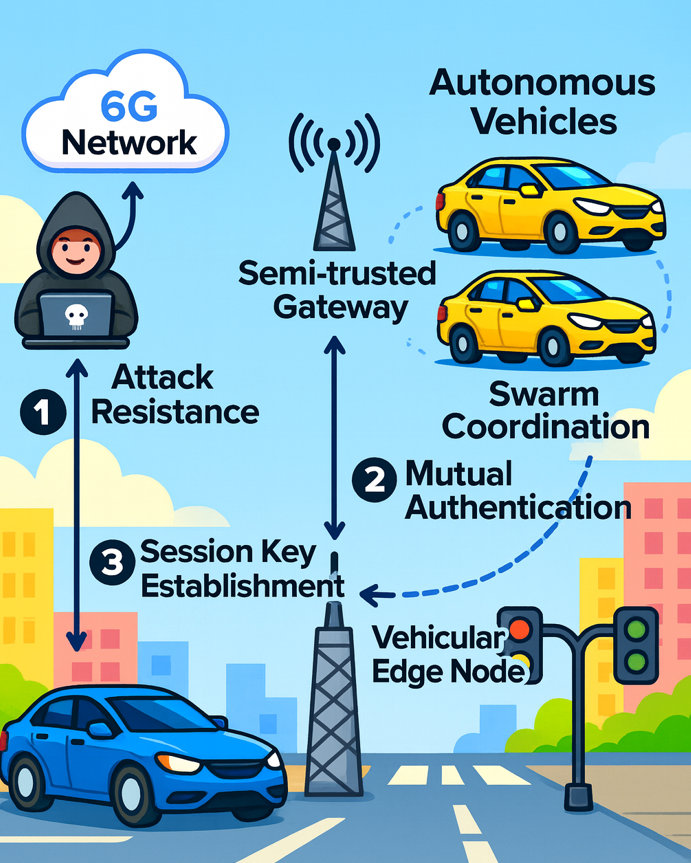}
\caption{QSAFE-V system model}
\label{fig:SM}
\end{figure}

AVs are intelligent vehicles capable of operating independently within a vehicular Tactile Internet environment. Each AV has various onboard sensors and communication units to support real-time decision-making and swarm coordination. Before participating in collaborative operations, AVs must authenticate with the VEN and obtain secure session credentials. Once authenticated, AVs can securely exchange tactile and haptic data with other swarm members, enabling cooperative driving tasks such as platooning, intersection negotiation, or obstacle avoidance.
VENs are edge computing units deployed near road infrastructure, such as at intersections or highways. They serve as the primary authentication and coordination points for nearby AVs. VENs perform context-aware behavioral analysis and manage lightweight reauthentication processes during high-mobility conditions. To ensure secure communication, VENs integrate QKD modules and act as quantum access points for the AVs.
GW is an intermediary device within the vehicular Tactile Internet, connecting the VENs and AVs to the central traffic management infrastructure. The GW facilitates secure key exchanges, relays entangled authentication tokens when needed, and coordinates swarm-level communication sessions. It is considered semi-trusted, with access limited to coordination logic but not complete session keys.
TA is a globally trusted entity that manages the long-term identities and credentials of AVs, VENs, and GWs. It distributes quantum-secured credentials, issues entangled token pairs for distributed authentication, and assists in bootstrapping the swarm coordination trust model.
%Tactile Internet  utilizes 6G-enabled wireless communication and tactile feedback channels to allow AVs to exchange real-time control, haptic, and sensory data. Vehicles interact in a coordinated and context-aware manner through low-latency quantum-secured wireless links, forming resilient vehicular swarms that operate autonomously without human intervention.

\subsection{Implicit Attack}
We studied the work in~\cite{Guo2019} and observed that many authentication protocols fail to account for the existence of implicit attacks during security analysis. This oversight leads to latent vulnerabilities in protocols that otherwise appear robust against traditional attacks. In other words, schemes may withstand explicit attacks such as replay or impersonation, but remain susceptible to combined or indirect forms of attack that exploit contextual and behavioral deviations in real-world environments.
In the context of the vehicular Tactile Internet, such implicit attacks can have catastrophic consequences—ranging from vehicle miscoordination to swarm control hijacking. Adversaries with multi-modal capabilities can leverage contextual shifts, vehicular mobility, or quantum computing advantages to exploit weak authentication models. For example, several classical schemes~\cite{Yadav2019, Zhang2020} were shown to lack session key secrecy and forward/backward security under implicit attack models. Additionally, some protocols fail to resist advanced impersonation attempts when adversaries capture vehicular context or inject adversarial data during vehicular handovers~\cite{He2018}.
\begin{table*}[ ]
\centering
\caption{Types of Known Attacks (Explicit Attacks)}
\begin{tabular}{|c|p{16.5cm}|}
\hline
\textbf{Item} & \textbf{Description} \\
\hline
KA1 & \textbf{DY Attacks.} The adversary $\mathcal{A}$ can eavesdrop, inject, or manipulate messages over public vehicular communication channels (e.g., V2V or V2I). \\
\hline
KA2 & \textbf{CK-I Attacks.} The adversary $\mathcal{A}$ gains temporary session-specific access to vehicular swarm credentials or tokens. \\
\hline
KA3 & \textbf{CK-II Attacks.} The adversary $\mathcal{A}$ has access to long-term identity keys or quantum token credentials stored on the vehicle. \\
\hline
KA4 & \textbf{Vehicular Node Attacks.} $\mathcal{A}$ retrieves sensitive swarm coordination data from VENs. (e.g., KA4-1: AV memory, KA4-2: VEN cache) \\
\hline
KA5 & \textbf{Insider Attacks.} Privileged internal adversaries (e.g., authorized AVs or infrastructure managers) can eavesdrop on or tamper with authentication messages. (KA5-1: vehicle registration phase, KA5-2: VEN storage) \\
\hline
KA6 & \textbf{Quantum-Aided Generic Attacks.} $\mathcal{A}$ performs offline key-guessing or token analysis using classical or quantum-assisted computation. \\
\hline
\end{tabular}
\label{table:explicit_attacks}
\end{table*}

Inspired by the findings in~\cite{Guo2019}, the QSAFE-V model incorporates a broader security evaluation strategy to address implicit attacks in the vehicular quantum Tactile Internet. In our approach, multiple security attributes are systematically analyzed in tandem with their corresponding implicit attack surfaces. For each security goal, we define an associated implicit attack that models indirect or behavioral threat vectors. If a protocol can resist all implicit attacks mapped to its claimed security goals, it is deemed robust against the strongest form of contextual compromise.
To this end, we analyze QSAFE-V against known attack frameworks, including the Dolev-Yao (DY)~\cite{Dolev1983} and Canetti-Krawczyk (CK)~\cite{Canetti2001} models, extended with context-aware adversarial capabilities. Furthermore, post-quantum attack surfaces such as quantum key tampering, quantum side-channel attacks, and entanglement hijacking are included. %These reflect real-world attacker models with access to quantum processors or vehicular behavioral data, thus providing a more realistic assessment of vehicular authentication in 6G tactile networks.
To facilitate secure protocol design, we define known attacks (explicit) and implicit attacks in Table~\ref{table:explicit_attacks} and Table~\ref{table:implicit_attacks}, respectively. These are mapped to nine security goals (SG1–SG9), covering both classical and quantum adversarial models.
\begin{table*}[ ]
\centering
\caption{Security Goals and Implicit Attacks}
\begin{tabular}{|c|p{5.5cm}|p{7.5cm}|}
\hline
\textbf{Item} & \textbf{Security Goals} & \textbf{Implicit Attacks} \\
\hline
SG1 & Mutual authentication and quantum-secure key & [KA1, KA2, KA5, KA6] or [KA1, KA3, KA5, KA6] \\
\hline
SG2 & Resistance to replay attacks & [KA1, KA2, KA5, KA6] or [KA1, KA3, KA5, KA6] \\
\hline
SG3 & Resistance to man-in-the-middle attacks & [KA1, KA2, KA5, KA6] or [KA1, KA3, KA5, KA6] \\
\hline
SG4 & Resistance to impersonation attacks & [KA1, KA2, KA4/KA5-2, KA6] or [KA1, KA3, KA4, KA6] \\
\hline
SG5 & Resistance to offline swarm credential guessing & [KA1, KA4-1, KA5, KA6] \\
\hline
SG6 & Session key confidentiality & [KA1, KA2, KA4/KA5-2, KA6] or [KA1, KA3, KA4, KA6] \\
\hline
SG7 & Perfect forward/backward secrecy & [KA1, KA3] \\
\hline
SG8 & Anonymity of vehicular identities & [KA1, KA2, KA3, KA5, KA6] \\
\hline
SG9 & Resistance to desynchronized attacks & [KA1, KA2, KA5, KA6] or [KA1, KA3, KA5, KA6] \\
\hline
\end{tabular}
\label{table:implicit_attacks}
\end{table*}

\subsection{Authentication Scheme}
In this subsection, we present a secure quantum-enhanced authentication framework designed for the wireless vehicular Tactile Internet, called QSAFE-V. The proposed QSAFE-V protocol is tailored for AV swarm coordination, operating under the edge-assisted Tactile Internet paradigm. QSAFE-V comprises three core phases: initialization, registration, and quantum-authenticated login and key agreement.
When the vehicular Tactile Internet system is deployed, initialization is first carried out by a fully TA, such as a central vehicular identity controller. All participating entities—including AVs, VENs, and roadside gateways (RGs)—must be securely registered with the TA prior to participation. %This process ensures that every vehicle receives a unique QIT, derived using QKD protocols or pre-shared entangled state resources.
AVs must complete a secure registration process with the TA before entering the vehicular swarm environment for the first time. The TA also maintains a registry of VENs and RGs responsible for relaying authentication signals and assisting in key agreement negotiations.
The initialization phase and the registration phase are both executed in physically secure or quantum-secured channels, possibly using satellite uplink or QKD optical fibers. Once an AV logs in through the login phase, it becomes a legitimate participant in the vehicular coordination framework.

During the authentication and key agreement phase, the AV and the target VEN (or neighboring AV) engage in mutual authentication. This phase leverages post-quantum cryptographic functions (e.g., Lattice-based signatures) and QITs, ensuring that the authentication is resistant to impersonation, forward/backward compromise, and quantum-assisted inference. Edge nodes and RGs facilitate the distribution of temporary session credentials without having direct access to long-term secrets.
All entities involved in QSAFE-V and their symbolic notations are described in Table~\ref{table:notations}.
\begin{table*}[ ]
\centering
\caption{Notations and Descriptions}
\label{table:notations}
\begin{tabular}{|c|p{11cm}|}
\hline
\textbf{Notation} & \textbf{Description} \\
\hline
$\texttt{QIT}_{AV}, \texttt{QIT}_{\texttt{VEN}}, \texttt{QIT}_{RG}$ & Quantum Identity Tokens for AVs, VENs, and RGs \\
\hline
$TA, \texttt{AV}_i, \texttt{VEN}_j, \texttt{RG}_k$ & Trusted Authority, Autonomous Vehicle, Vehicular Edge Node, Roadside Gateway \\
\hline
$\texttt{RID}_{AV}, \texttt{RID}_{\texttt{VEN}}, \texttt{RID}_{RG}$ & Real identity of AV, VEN, and RG \\
\hline
$\texttt{PID}_{AV}, \texttt{PID}_{\texttt{VEN}}, \texttt{PID}_{RG}$ & Pseudo-identity for anonymization and unlinkability \\
\hline
$\texttt{PW}_i, H\texttt{PW}_i$ & Password and hashed password for AV\textsubscript{i} \\
\hline
$\texttt{Enc}_{PK}()$ & Asymmetric encryption and decryption using public and private keys (post-quantum primitives) \\
\hline
$n_i, r_i, \tau_{RG}, r_j$ & Random numbers and fresh nonces \\
\hline
$T_i, \Delta T$ & Timestamp and maximum allowable transmission delay \\
\hline
$SK_{TA}, sk_i, sk_j$ & Secret keys of TA and the AVs/VENs \\
\hline
$\texttt{PK}_{AV}, \texttt{PK}_{\texttt{VEN}}, \texttt{PK}_{RG}$ & Public keys of AVs, VENs, and RGs \\
\hline
$SSK_{i,j}, SK_{ij}$ & Shared session key between communicating parties \\
\hline
$\texttt{SKV}_i$ & Session key validator (used for key confirmation and verification) \\
\hline
$h()$ & Cryptographic hash function (quantum-resistant, e.g., SHA3 or SPHINCS+) \\
\hline
$\parallel, \oplus$ & Concatenation and XOR operations \\
\hline
\end{tabular}
\end{table*}

\subsubsection{Initialization Phase}
During this phase, the TA performs the initialization of the vehicular Tactile Internet system for quantum-secured autonomous swarm communication.

\textbf{Step I1:} TA begins by generating a post-quantum cryptographic (PQC) key pair using a secure lattice-based scheme such as CRYSTALS-Kyber or Dilithium. The TA selects a sufficiently strong private key $sk_{TA}$ and derives the corresponding public key $\texttt{PK}_{TA}$. These keys are used for initial encryption and identity attestation across the swarm domain.

\textbf{Step I2:} TA defines a secure QIT issuance process for all participating AVs, VENs, and RGs. Each QIT is bound to the real identity of the device and embedded using either:
- pre-shared entangled photon pairs (in QKD networks), or
- classical quantum-safe tokens generated via secure hash-then-sign primitives.

\textbf{Step I3:} TA selects a quantum-resistant hash function $h()$ from the SHA-3 or SPHINCS+ family to ensure resilience against pre-image and collision attacks. This function is used in generating pseudo-identities and session key validators.
Finally, TA stores the master private key $sk_{TA}$ in a secure hardware enclave and publishes the system-wide public parameters as:
$\{ \texttt{PK}_{TA}, h(), \texttt{QIT}_{format}, QSAFE\}$
These parameters are broadcast to all legitimate entities within the QSAFE-V network, enabling downstream registration and authenticated swarm coordination.

\subsubsection{Registration Phase}
This phase is conducted offline by the TA to register all participating entities: ($\texttt{RG}_k$, $k = 1, 2, ..., n_{RG}$), autonomous vehicles ($\texttt{AV}_i$, $i = 1, 2, ..., n_{AV}$), and vehicular edge nodes ($\texttt{VEN}_j$, $j = 1, 2, ..., n_{\texttt{VEN}}$). The following subsections detail the registration procedure for each entity.

\textit{1) RG Registration}

\textbf{Step RG1:} TA selects the real identity $\texttt{RID}_{RG}$ for the RG, computes its pseudo-identity $\texttt{PID}_{RG} = h(\texttt{RID}_{RG} \parallel \kappa)$ using a secure post-quantum hash, and generates a private key $sk_{RG}$ using a lattice-based key generation scheme.

\textbf{Step RG2:} TA issues the registration certificate as $\{\texttt{PID}_{RG}, sk_{RG}\}$ and stores it securely at the RG for future authentication operations.

\textit{2) Autonomous Vehicle Registration}

\textbf{Step AV1:} AV$_i$ selects its real identity $\texttt{RID}_{\texttt{AV}_i}$ and local password $\texttt{PW}_i$, computes the encrypted registration request $\texttt{Req}_i = \texttt{Enc}_{\texttt{PK}_{TA}}(\texttt{RID}_{\texttt{AV}_i})$, and sends it to TA via a protected control channel.

\textbf{Step AV2:} Upon receipt, TA decrypts $\texttt{Req}_i$ and verifies $\texttt{RID}_{\texttt{AV}_i}$. It then assigns a temporary quantum identity token $\texttt{QIT}_i$ and generates a lattice-based key pair $\{sk_i, \texttt{PK}_i\}$ for AV$_i$. The TA then issues a pseudo-identity $\texttt{PID}_{\texttt{AV}_i} = h(\texttt{RID}_{\texttt{AV}_i} \parallel \kappa)$ and sends $\{\texttt{QIT}_i, \texttt{PK}_i, \texttt{PID}_{RG}\}$ back to AV$_i$.

\textbf{Step AV3:} AV$_i$ generates a random nonce $n_i \in \mathbb{Z}_q^*$ and stores the registration parameters $\{\texttt{QIT}_i, n_i, sk_i, \texttt{PID}_{RG} \oplus \texttt{PID}_{\texttt{AV}_i}\}$ securely within its onboard cryptographic module. The public key $\texttt{PK}_i$ is exposed for mutual authentication.

\textit{3) Vehicular Edge Node Registration}

\textbf{Step VEN1:} For each deployed $\texttt{VEN}_j$, TA assigns a unique identity $\texttt{RID}_{\texttt{VEN}_j}$, calculates the pseudo-identity $\texttt{PID}_{\texttt{VEN}_j} = h(\texttt{RID}_{\texttt{VEN}_j} \parallel \kappa)$, and generates a challenge nonce $C_{\texttt{VEN}_j}$.

\textbf{Step VEN2:} TA computes a post-quantum key pair $\{sk_j, \texttt{PK}_j\}$ and binds it to $\texttt{VEN}_j$.

\textbf{Step VEN3:} TA finalizes registration by issuing credentials $\{\texttt{PID}_{\texttt{VEN}_j}, sk_j, h(\texttt{PID}_{GW}), C_{\texttt{VEN}_j}\}$ to $\texttt{VEN}_j$ before deployment.

\begin{algorithm}[ ]
\caption{Registration Phase of AV}
\KwIn{$\texttt{RID}_i$, $\texttt{PW}_i$, \newline $MR_2 = \{\texttt{QIT}_i, sk_i, \texttt{PK}_i, \texttt{PID}_{RG} \}$}
\KwOut{$MR_1 = \texttt{Req}_i$ or \textbf{false}}
\Begin{
    $\texttt{AV}_i$ inputs $\texttt{RID}_i$ and $\texttt{PW}_i$, and computes registration request \\
    \Indp $\texttt{Req}_i = \texttt{Enc}_{\texttt{PK}_{TA}}(\texttt{RID}_i)$ \;
    \Indm Send $MR_1 = \texttt{Req}_i$ to TA \;
    
    TA decrypts $\texttt{Req}_i$ to get $\texttt{RID}_i$ \;
    TA generates unique quantum identity token $\texttt{QIT}_i$ \;
    Generate key pair: $sk_i \in \mathbb{Z}_q^*$, $\texttt{PK}_i = f(sk_i)$ \;
    Compute $\texttt{PID}_i = h(\texttt{RID}_i \parallel \kappa)$ \;
    Send $MR_2 = \{\texttt{QIT}_i, sk_i, \texttt{PK}_i, \texttt{PID}_{RG}\}$ to $\texttt{AV}_i$ \;

    \eIf{$\texttt{AV}_i$ receives $MR_2$}{
        Generate random nonce $n_i \in \mathbb{Z}_q^*$ \;
        Compute: \\
        \Indp $A = h(\texttt{RID}_i \parallel \texttt{PW}_i) \oplus n_i$ \;
        $\texttt{PID}_i = h(\texttt{RID}_i \parallel n_i)$ \;
        \Indm Store $\{A, \texttt{QIT}_i, n_i, sk_i, \texttt{PID}_{RG} \oplus \texttt{PID}_i\}$ in secure module \;
    }{
        \Return \textbf{false} \;
    }
}
\end{algorithm}

\subsubsection{Autonomous Vehicle Login and Identity Authentication Key Agreement Phase}
In this phase, a registered autonomous vehicle $\texttt{AV}_i$ attempts to establish a secure session key with a $\texttt{VEN}_j$, assisted by ($\texttt{RG}_k$), after a successful login. The session key obtained is used to enable authenticated swarm coordination and URLLC. The protocol progresses through the following steps:

\textbf{Step L1:} When $\texttt{AV}_i$ enters its credentials $\{\texttt{RID}_i, \texttt{PW}_i\}$ into its onboard terminal, it computes $n_i' = A \oplus h(\texttt{RID}_i \parallel \texttt{PW}_i)$ and verifies its stored registration value. If validation is successful, the login proceeds.

\textbf{Step QRM1:} $\texttt{AV}_i$ generates a random nonce $r_i \in \mathbb{Z}_q^*$ and current timestamp $T_1$, calculates its pseudo-identity $\texttt{PID}_i = h(\texttt{RID}_i \parallel n_i)$ and its quantum identity token $\texttt{QIT}_i$. Then, it computes:
$S_i = \texttt{PQK}_i,\quad R_i = \texttt{Enc}_{\texttt{PK}_{\texttt{VEN}_j}}(\texttt{QIT}_i \parallel \texttt{PID}_i \parallel T_1)$
The message $M_1 = \{r_i, R_i, \texttt{PK}_i, T_1\}$ is transmitted to the $\texttt{RG}_k$ via an authenticated channel.

\textbf{Step QRM2:} Upon receiving $M_1$ at time $T_1'$, the $\texttt{RG}_k$ first checks time validity $|T_1' - T_1| < \Delta T$. If satisfied, it verifies $\texttt{PK}_i$ and decrypts $R_i$ to retrieve $\texttt{QIT}_i$ and $\texttt{PID}_i$. Then, $\texttt{RG}_k$ generates a nonce $r_{RG}$ and timestamp $T_2$, computes:
$D = r_{RG} \oplus h(\texttt{PID}_{RG} \parallel T_2), \quad C = \texttt{Enc}_{\texttt{PK}_{\texttt{VEN}_j}}(D \parallel T_2)$
and sends $M_2 = \{r_i, C, T_2\}$ to $\texttt{VEN}_j$.

\textbf{Step QRM3:} $\texttt{VEN}_j$ validates timestamp $T_2$ and decrypts $C$ to get $D$. It then verifies:
$\texttt{PID}_{RG} = D \oplus h(r_i \parallel \texttt{PK}_i)$
If valid, $\texttt{VEN}_j$ generates new nonce $r_j$ and timestamp $T_3$, computes:
$F = r_j \oplus h(\texttt{PID}_i \parallel T_3)$
and sends $M_3 = \{r_j, \texttt{SKV}_i, F, T_3\}$ to the $\texttt{RG}_k$.

\textbf{Step QRM4:} $\texttt{RG}_k$ receives $M_3$ and forwards it to $\texttt{AV}_i$ at time $T_4$. The vehicle checks $|T_4 - T_3| < \Delta T$. Then it calculates:
$F' = r_j \oplus h(\texttt{PID}_i \parallel T_3)$
If verification passes, $\texttt{AV}_i$ calculates the session key:
$SK_{ij} = h(S_i \oplus r_j \parallel T_3), \quad \texttt{SKV}_i = h(SK_{ij} \parallel T_3)$

\textbf{Step QRM5:} Finally, $\texttt{AV}_i$ updates its temporary identity:
$\texttt{QIT}_i' = h(\texttt{QIT}_i \parallel \texttt{SKV}_i \parallel T_4)$
and confirms successful mutual authentication.

Through these steps, $\texttt{AV}_i$ and $\texttt{VEN}_j$ securely establish a shared session key $SK_{ij}$ via the $\texttt{RG}_k$, ensuring quantum-safe swarm interaction.
\begin{algorithm}[ ]
\caption{Quantum-Secured Authentication}
\KwIn{$ID'_i$, $PW'_i$, Quantum-Enhanced Message $QM_4 = \{R_j, \texttt{SKV}_j, \mathcal{Q}, C, F, T_3, T_4\}$ from $GW$}
\KwOut{$QM_1 = \{R_i, \mathcal{B}, Pub_i, T_1\}$ or \textbf{fail}}

\textbf{Step Q1 – Quantum Login:}

\Begin{
  Calculate $n'_i = A \oplus h(ID'_i \| PW'_i)$\;
  \If{$n'_i == n_i$}{
    Generate random $r_i \in \mathbb{Z}_q^*$ and current timestamp $T_1$\;
    Calculate:
    \begin{itemize}
        \item $\texttt{PID}_i = h(\texttt{TID}_i \| n_i)$
        \item $S_i = h(r_i \| \texttt{TID}_i \| \texttt{PID}_i \| T_1)$
        \item $R_i = S_i \cdot P$
        \item $\mathcal{B} = S_i + k_i$
    \end{itemize}
    Send $QM_1 = \{R_i, \mathcal{B}, Pub_i, T_1\}$ to $GW$\;
  }
  \Else{return \textbf{fail}}
}

\vspace{0.2cm}
\textbf{Step Q5 – Quantum Key Validation (Wait for $GW$ response):}

\Begin{
Receive $QM_4 = \{R_j, \texttt{SKV}_j, \mathcal{Q}, C, F, T_3, T_4\}$ from $GW$\;

\If{$|T_4' - T_4| < \Delta T$}{
  Calculate $\texttt{PID}_{GW} = \texttt{PID}_i \oplus (\texttt{PID}_{GW} \oplus \texttt{PID}_i)$\;
  $V^{**}_{GW} = \mathcal{Q} \oplus h(R_i \| Pub_j)$\;

  \If{$h(\texttt{PID}_{GW} \| T_4) == F \oplus V^{**}_{GW}$}{
    $SK_i = h(S_i \cdot R_j \| V^{**}_{GW} \| T_3)$\;
    $\texttt{SKV}_i = h(SK_i \| T_3)$\;
    \If{$\texttt{SKV}_i == \texttt{SKV}_j$}{
      Update $\texttt{TID}_i = h(\texttt{TID}_i \| SK_i \| T_4)$\;
    }
    \Else{return \textbf{fail}}
  }
  \Else{return \textbf{fail}}
}
\Else{return \textbf{fail}}
}
\end{algorithm}

\section{Security Analysis}

According to \cite{mohsin2016iotsat}, combining formal security proofs with non-formal security analysis is essential when evaluating the robustness of authentication frameworks in next-generation communication systems. Classical security proofs often fail to comprehensively capture advanced adversarial capabilities—especially in quantum-driven networks.
To address this gap, we introduce a dual-layer analysis approach for the proposed QSAFE-V framework. First, we redefine and expand the concept of implicit quantum attacks, which exploit the interdependencies of QITs, session randomness, and entanglement leakage. This extension is necessary to ensure provable resilience in the presence of quantum adversaries equipped with both classical and quantum eavesdropping capabilities.

Building upon this foundation, we extend the widely accepted ROR oracle model~\cite{becerra2017relation} to support quantum-side oracle access and simulate adversary behavior in hybrid quantum-classical environments. The resulting proof validates that QSAFE-V ensures indistinguishability of session keys under chosen-session and adaptive quantum attacks.
Furthermore, we demonstrate that QSAFE-V effectively mitigates both implicit and explicit threats, including quantum impersonation, session key inference, and entanglement hijacking. These evaluations confirm QSAFE-V's ability to achieve the listed security goals even under adversarial quantum computation.
Simulation experiments are conducted using the AVISPA tool (Automated Validation of Internet Security Protocols and Applications) with quantum-aware extensions, further affirming the correctness and robustness of QSAFE-V under symbolic and protocol-level attacks.

\subsection{Extended ROR Model}

Formal security analysis based on the ROR model is a powerful method to prove the security of cryptographic authentication schemes. To account for the quantum-enhanced threat space and the concept of implicit quantum attacks, we extend the classical ROR model to suit the QSAFE-V protocol context. The components of the extended model are as follows:
Participants in QSAFE-V include the autonomous vehicle $(\texttt{AV}_i)$, $(\texttt{VEN}_i)$,  $(\texttt{GW}_i)$, and $(TA)$. Protocol instances are denoted respectively as $\Pi^{u}_{AV}$, $\Pi^{v}_{VEN}$, and $\Pi^{w}_{GW}$, where $u, v, w$ index the instances.
 An instance $\Pi^{x}$ transitions to the accepted mode after completing all message exchanges and validating the QIT. A unique session identifier $(\texttt{sid})$ is derived by concatenating all exchanged messages and timestamps.
 Two protocol instances $\Pi^{x_1}_i$ and $\Pi^{x_2}_j$ are considered partnered if:
both are in the accepted mode; they authenticate each other and derive identical session keys; and they share the same session identifier $(\texttt{sid})$.
A session is fresh if the session key has not been exposed via the \texttt{Reveal} query or through any quantum side-channel leakage vector.
 The following oracle queries model the power of a quantum-capable adversary $\mathcal{A}$:

 \texttt{ExecuteKA1}$(\Pi^{u}_{AV}, \Pi^{v}_{VEN}, \Pi^{w}_{GW})$: $\mathcal{A}$ eavesdrops on the public exchange between AV, VEN, and GW.

 \texttt{SendKA1}$(\Pi^{x}, m)$: $\mathcal{A}$ can send, replay, or modify a message $m$ to the instance $\Pi^x$ and observe the output.

\texttt{CorruptKA2}$(\Pi^{x})$: Leaks temporary session data (e.g., random values, timestamps, ephemeral QIT hashes).

 \texttt{CorruptKA3}$(\Pi^{x})$: Reveals long-term secrets (e.g., lattice-based private keys, pre-shared QKD credentials).

 \texttt{CorruptKA4\_1}$(\Pi^{u}_{AV})$: Reveals stored QITs and device credentials of the AV.

 \texttt{CorruptKA4\_2}$(\Pi^{v}_{VEN})$: Reveals entanglement token mappings or encoded quantum ID signatures at the VEN.

 \texttt{CorruptKA5\_1}$(\Pi^{u}_{AV})$: Eavesdrops on AV's registration phase with TA.

 \texttt{CorruptKA5\_2}$(\Pi^{w}_{GW})$: Reveals pre-configured security credentials stored in the gateway.

 \texttt{Reveal}$(\Pi^x)$: Leaks the session key established by $\Pi^x$ and its partner.

\texttt{Test}$(\Pi^x)$: Models session key indistinguishability. The challenger flips a hidden bit $c \in \{0,1\}$ and returns:
    \begin{itemize}
        \item real session key if $c=1$ and the session is fresh;
        \item random key otherwise.
    \end{itemize}
    $\mathcal{A}$ must guess $c$ with success probability non-negligibly better than $1/2$.

\begin{figure}[!ht]
\centering
\includegraphics[width=2.9in]{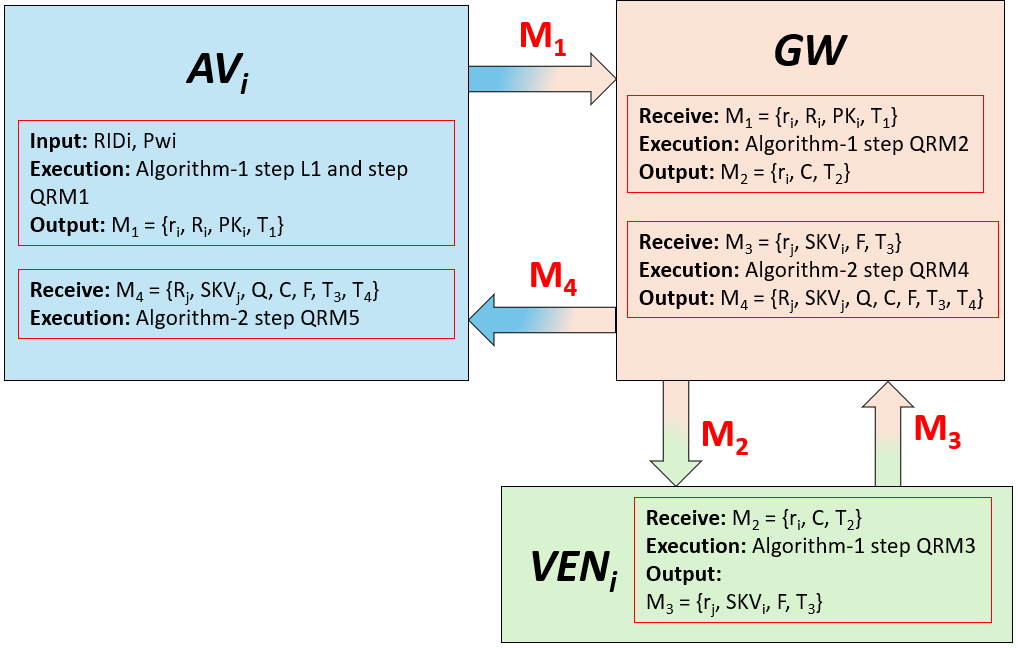}
\caption{Flowchart of the algorithm.}
\label{fig:flowchart}
\end{figure}
\subsection{Semantic Security of the Session Key}
In the extended ROR model, the adversary $\mathcal{A}$ is required to distinguish between the real session key of a returned instance and a randomly generated key of equal length. To achieve this, $\mathcal{A}$ is allowed to issue a series of queries including \texttt{Execute}, \texttt{Send}, \texttt{Reveal}, \texttt{Corrupt}, and \texttt{Test} multiple times within the defined security game.
At the conclusion of the game, $\mathcal{A}$ outputs a guess $c'$ for the challenge bit $c$. The adversary is considered successful if $c' = c$. Let $\texttt{Succ}$ denote the event where $\mathcal{A}$ correctly guesses $c$. Then, the advantage of $\mathcal{A}$ in breaking the semantic security of the session key for authentication scheme $\mathcal{S}$ is defined as:
$\texttt{Adv}_{\mathcal{S}}(\mathcal{A}) = \left| 2 \cdot \Pr[\texttt{Succ}] - 1 \right|$.
The authentication scheme $\mathcal{S}$ is said to achieve semantic security under the extended ROR model if for any probabilistic polynomial-time adversary $\mathcal{A}$, the advantage satisfies $\texttt{Adv}_{\mathcal{S}}(\mathcal{A}) \leq \epsilon$, where $\epsilon$ is a negligible value that converges to zero as a function of the security parameter.
In the QSAFE-V protocol, all participants, including the adversary $\mathcal{A}$, have access to a collision-resistant one-way hash function $h(\cdot)$ modeled as a random oracle. This ensures that preimage and collision attacks remain computationally infeasible even under quantum constraints.

\subsection{Formal Security Proof Using the Extended ROR Model}

In this section, we demonstrate the semantic security of the proposed QSAFE-V protocol. Under the extended ROR model, the maximum advantage of an adversary $\mathcal{A}$ in compromising the session key security of QSAFE-V within probabilistic polynomial time is bounded as follows:

\begin{equation}
\texttt{Adv}^{max}_{\mathcal{S}}(\mathcal{A}) \leq \frac{q^2_h}{|\texttt{Hash}|} + 2\texttt{Adv}^{\mathcal{LWE}}_{\mathcal{S}}(\mathcal{A}),
\end{equation}

where $q_h$ denotes the number of hash oracle queries, $|\texttt{Hash}|$ is the range of the collision-resistant hash function modeled as a random oracle, and $\texttt{Adv}^{\mathcal{LWE}}_{\mathcal{S}}(\mathcal{A})$ is the advantage of $\mathcal{A}$ in solving the LWE problem, which underpins the lattice-based primitives in QSAFE-V.

\textbf{Proof:} We define a sequence of games $G_i$ $(i = 0, 1, \dots, 6)$ to analyze the advantage of $\mathcal{A}$ in distinguishing a real session key from a random one. Let $\texttt{Succ}_i$ denote the success probability of $\mathcal{A}$ in game $G_i$.

\textbf{Game $G_0$:} This game simulates a real attack scenario where $\mathcal{A}$ has no prior knowledge. The advantage is:

\begin{equation}
\label{Eq_1}
\texttt{Adv}_{\mathcal{S}}(\mathcal{A}) = |2 \cdot Pr[\texttt{Succ}_0] - 1|.
\end{equation}

\textbf{Game $G_1$:} $\mathcal{A}$ executes \texttt{ExecuteKA1} and \texttt{SendKA1} to intercept messages between AV$_i$, VEN$_j$, and GW$_k$:
\[
M_1 = \{\texttt{AV}_i, \texttt{Pub}_i, \texttt{QIT}_i\}, \quad M_2 = \{\texttt{VEN}_i, C, D, \texttt{QIT}_j\},
\]
\[
M_3 = \{\texttt{GW}_i, \texttt{SKV}_{i,j}, T_3, T_4\}.
\]
Then $\mathcal{A}$ invokes \texttt{CorruptKA2} to obtain session-specific data.

\textbf{Game $G_2$:} $\mathcal{A}$ lacks knowledge of secret values like $\texttt{PID}_{GW}$, preventing accurate computation of the session key:
\begin{equation}
\label{Eq_2}
Pr[\texttt{Succ}_2] - Pr[\texttt{Succ}_1] = 0.
\end{equation}

\textbf{Game $G_3$:} $\mathcal{A}$ invokes \texttt{CorruptKA3} to obtain long-term lattice-based private keys $k_i$, $k_j$, and $k_{GW}$.

\textbf{Game $G_4$:} $\mathcal{A}$ uses \texttt{CorruptKA4} to obtain stored parameters from AV or VEN. Still cannot compute shared QIT state correctly due to the entanglement dependency.

\textbf{Game $G_5$:} Attacker uses \texttt{CorruptKA5} to retrieve registration phase metadata. Since registration does not reveal QIT hashes, we have:

\begin{equation}
\label{Eq_3}
Pr[\texttt{Succ}_5] - Pr[\texttt{Succ}_4] = 0.
\end{equation}

\textbf{Game $G_6$:} $\mathcal{A}$ attempts to guess bit $c$ in \texttt{Test} query without full knowledge of lattice-based session key derivation:

\begin{equation}
\label{Eq_4}
Pr[\texttt{Succ}_6] = \frac{1}{2}.
\end{equation}

From (\ref{Eq_1}) through (\ref{Eq_4}), we conclude:

\begin{equation}
\texttt{Adv}_{\mathcal{S}}(\mathcal{A}) \leq \frac{q^2_h}{|\texttt{Hash}|} + 2\texttt{Adv}^{\mathcal{LWE}}_{\mathcal{S}}(\mathcal{A}).
\end{equation}

\textbf{Implicit Attack Paths.}
\begin{itemize}
  \item Implicit attack [KA1, KA2, KA5-2, KA6]: $\mathcal{A}$ fails due to lack of $\texttt{Reg}_{GW}$ or QIT seed.
  \item Implicit attack [KA1, KA3, KA6]: Even with lattice private keys, $\mathcal{A}$ lacks entropy alignment for QIT token.
  \item Implicit attack [KA1, KA4, KA6]: Without synchronized $\texttt{TID}_i$ and pre-shared entangled QIT, session computation fails.
\end{itemize}

In all paths, $\mathcal{A}$ must break LWE or hash collisions to succeed. Therefore, we conclude:

\begin{equation}
\texttt{Adv}^{max}_{\mathcal{S}}(\mathcal{A}) \leq \frac{q^2_h}{|\texttt{Hash}|} + 2\texttt{Adv}^{\mathcal{LWE}}_{\mathcal{S}}(\mathcal{A}).
\end{equation}

\subsection{Enhanced Non-Formal Security Analysis}
This section combines the mapping between security attributes and implicit attacks listed in Table~\ref{table:implicit_attacks}, using an enhanced non-formal security analysis methodology to demonstrate that QSAFE-V satisfies all defined security properties—proving its resilience against quantum-driven implicit attacks.

\subsubsection{Mutual Authentication and Key Establishment (SG1)}
SG1 ensures mutual authentication and session key establishment among $\texttt{AV}_i$, $\texttt{VEN}_i$, and $\texttt{GW}_i$ under implicit attacks [KA1, KA2, KA5, KA6] and [KA1, KA3, KA5, KA6].
If an adversary $\mathcal{A}$ can intercept:
$M_1 = \{\texttt{AV}_i, \texttt{Pub}_i, \texttt{QIT}_i\}, \quad M_2 = \{\texttt{VEN}_i, C, D, T_2\}, \quad M_3 = \{\texttt{GW}_i, \texttt{SKV}_{i,j}, T_3, T_4\}$,
they gain access to public values and ephemeral session parameters. However, $\mathcal{A}$ must also forge or modify the quantum identity token $\texttt{QIT}_i$ and secret values like $\texttt{PID}_{GW}$ and $\texttt{TID}_i$ to construct a valid session.
Even under attacks [KA1, KA3, KA5, KA6], the session key:
$SK_i = h(k_i(\tau_i \| \texttt{TID}_i \| \texttt{PID}_i \| T_1)) \cdot P \quad \text{or} \quad SK_j = h(k_j(\texttt{QIT}_j \| \texttt{PID}_{GW} \| T_3))$,
remains secure due to the adversary's inability to compute missing private parameters. Thus, SG1 is satisfied.

\subsubsection{ Resistance to Replay Attacks (SG2)}
SG2 addresses replay resilience under [KA1, KA2, KA5, KA6] and [KA1, KA3, KA5, KA6]. Dynamic timestamps $T_1$ and session-specific QIT tokens ensure uniqueness. Any attempt by $\mathcal{A}$ to replay $M_1$ or $M_2$ will be rejected unless all components (e.g., $\texttt{PID}_i$, $\texttt{TID}_i$) match and remain fresh, which $\mathcal{A}$ cannot guarantee.

\subsubsection{Resistance to Man-in-the-Middle Attacks (SG3)}
SG3 is fulfilled by using authenticated quantum channels and session-specific parameters (e.g., $\texttt{QIT}_i$, $\texttt{TID}_i$). Adversaries under [KA1, KA2, KA5, KA6] cannot modify or forge exchanged messages due to their dependency on quantum authentication hashes and randomness tied to the LWE key material.

\subsubsection{Resistance to Impersonation Attacks (SG4)}
Adversaries [KA1, KA3/KA4/KA5, KA6] may extract some long-term keys but lack $\texttt{PID}_{GW}$ and cannot reconstruct $\texttt{QIT}_j$ or valid entangled tokens. Impersonation toward $GW$ or $\texttt{VEN}_i$ fails due to QIT hash verifications and misalignment in entangled quantum token decoding.

\subsubsection{Resistance to Offline Password-Guessing Attacks (SG5)}
Attackers [KA1, KA4-1, KA6] may try to guess passwords using public messages and intercepted data, but QSAFE-V integrates one-time pads and entropy from QKD/PUF devices (e.g., $PUF_{\texttt{AV}_i}$) in session derivation. The hashes such as $h(ID_i \| r_i)$ and session randomness make password guesses indistinguishable from random noise.

\subsubsection{Session Key Security (SG6)}
Even under implicit attacks [KA1, KA3], the adversary cannot deduce session keys of previous or future sessions. This is due to the dynamic change of $T$, $r$, and $\texttt{QIT}$ values, as well as the entropy embedded in:
$SK = h(k(\tau \| \texttt{PID}_{GW} \| T)) - P \| k(\texttt{PID}_{GW}) \| TID.$
SG6 is therefore satisfied under quantum-resilient design.

\subsubsection{Perfect Forward/Backward Secrecy (SG7)}
Even if an adversary later learns long-term secrets, they cannot reconstruct past or future session keys due to their dependency on ephemeral values $\texttt{TID}_i$, $r_i$, and $PUF_{VEN}$. Implicit attacks [KA1, KA3] are insufficient without recomputing one-time parameters.

\subsubsection{Anonymity (SG8)}
Under [KA1, KA2, KA5, KA6], adversaries may know $M_1$, $M_2$, $r_i$, and $\texttt{TID}_i$, but cannot compute $\texttt{PID}_i = h(ID_i \| \tau_i)$ due to unknown $\tau_i$. Since $ID_i$ is only known to TA and AV$_i$, anonymity is preserved. Offline guesses fail due to LWE-based security and hash protection.

\subsubsection{Resistance to Desynchronized Attacks (SG9)}
Desynchronization-resistance in QSAFE-V is maintained by embedding $\texttt{TID}_i$ with quantum timestamping and secure storage mechanisms. Even if messages are delayed or intercepted, quantum verifiers check for entanglement validity and key freshness, which mitigates sync-based disruptions.

\begin{table}[!t]
\centering
\caption{Comparison of Performance}
\label{tab:performance}
\begin{tabular}{|c|cccccccc|c|c|c|c|}
\hline
 \textbf{Ref.} & \textbf{SG1} & \textbf{SG2} & \textbf{SG3} & \textbf{SG4} & \textbf{SG5} & \textbf{SG6} & \textbf{SG7} & \textbf{SG8}  \\
\hline
\cite{masud2021} & \checkmark & \checkmark & $\times$ & \checkmark & \checkmark & $\times$ & $\times$ & \checkmark\\
\cite{lee2020} & \checkmark & \checkmark & \checkmark & \checkmark & $\times$ & $\times$ & $\times$ & \checkmark\\
%\cite{kim2019} & \checkmark & \checkmark & \checkmark & \checkmark & $\times$ & \checkmark & \checkmark & $\times$ \\
%\cite{zhou2019} & \checkmark & \checkmark & \checkmark & \checkmark & \checkmark & $\times$ & $\times$ & $\times$\\
\cite{kaur2020} & \checkmark & \checkmark & \checkmark & \checkmark & $\times$ & $\times$ & $\times$ & $\times$\\
%\cite{hu2021} & \checkmark & $\times$ & $\times$ & \checkmark & \checkmark & \checkmark & \checkmark & $\times$\\
\cite{kamarudin2017development} & \checkmark & \checkmark & \checkmark & \checkmark & \checkmark & \checkmark & \checkmark & $\times$\\
\textbf{QSAFE-V} & \checkmark & \checkmark & \checkmark & \checkmark & \checkmark & \checkmark & \checkmark & \checkmark\\
\hline
\end{tabular}
\end{table}

\section{Results}
This section focuses on the performance evaluation result for our proposed QSAFE-V and the simulation of attack resistance.
We show our proposed QSAFE-V performance for all implicit attacks.
The overall comparison is shown in Table~\ref{tab:performance}.
As can be seen from Table~\ref{tab:performance}, lightweight protocols often face challenges in achieving multiple security goals under implicit attacks. None of the referenced lightweight schemes fully satisfy SG1, SG2, SG3, SG4, SG5, SG6, SG7, and SG8, which represent the fundamental security objectives of any authentication mechanism. Those Implicit attacks substantially undermine the robustness and reliability of secure communications, especially in critical domains such as autonomous vehicles and smart healthcare systems.
For schemes using public key cryptography, only the proposed QSAFE-V has achieved all the security goals in the table under implicit attacks and ensured both quantum-resilient authentication and low overhead. 
\begin{figure}[H]
\centering
\includegraphics[width=0.80\linewidth]{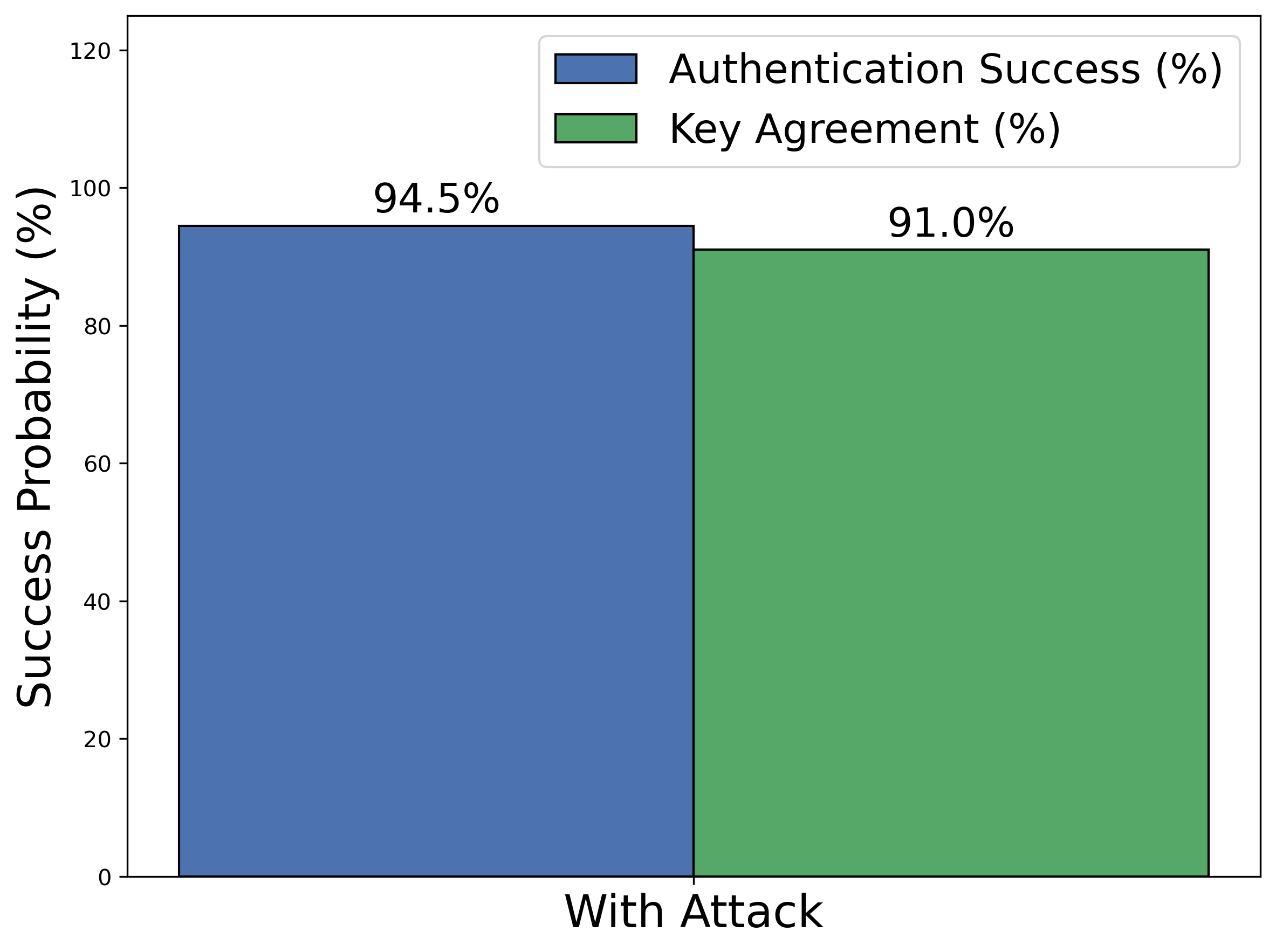}
\caption{Mutual Authentication and Key establishment for QSAFE-V.}
\label{fig:mutual_authentication}
\end{figure}

\subsection{Experimental Simulation Analysis Based on QISKIT.}
We utilize IBM-Qiskit to simulate eavesdropping, man-in-the-middle (MITM), and other implicit attacks on quantum communication channels to verify the robustness and key confidentiality of the proposed QSAFE-V scheme with URLLC in the Tactile Internet conditions. The simulation results are shown below.
Fig.~\ref{fig:mutual_authentication} shows the performance of the proposed mutual authentication and key agreement scheme under all attack conditions. Two bars represent the authentication success rate and key agreement success rate when an adversary attempts to disrupt or impersonate legitimate sessions. This validates security goal (SG1), ensuring that legitimate entities (AVi, VENi, GWi) can mutually authenticate and derive a consistent session key despite implicit attacks. The attacker cannot reproduce the QITi, preserving session confidentiality and authenticity.

\begin{figure}[H]
     \centering
     \begin{subfigure}[b]{0.45\textwidth}
         \centering 
         \includegraphics[width=1\textwidth]{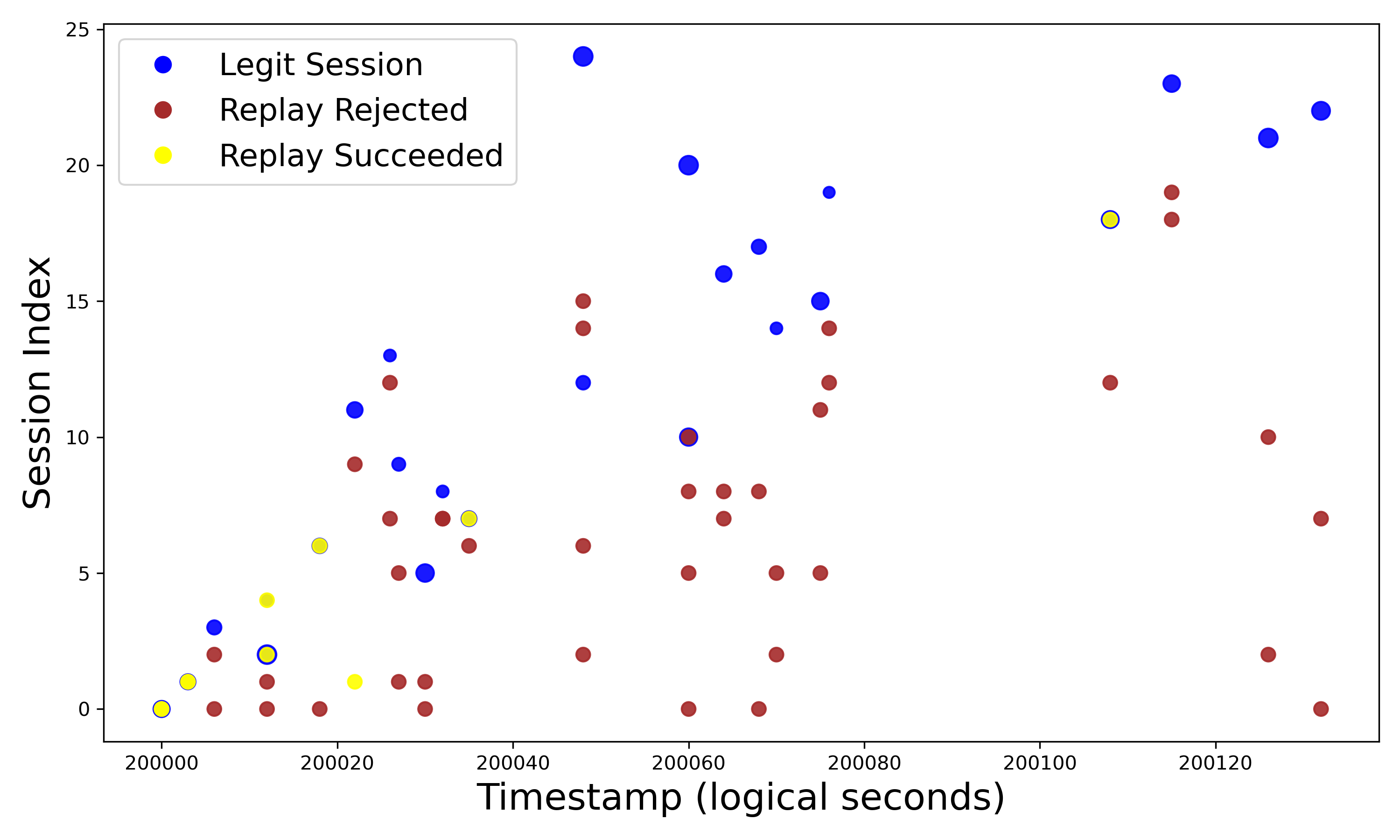}
         \caption{}
         \label{fig:..}
     \end{subfigure}
     \hfill
     \begin{subfigure}[b]{0.45\textwidth}
         \centering
         \includegraphics[width=\textwidth]{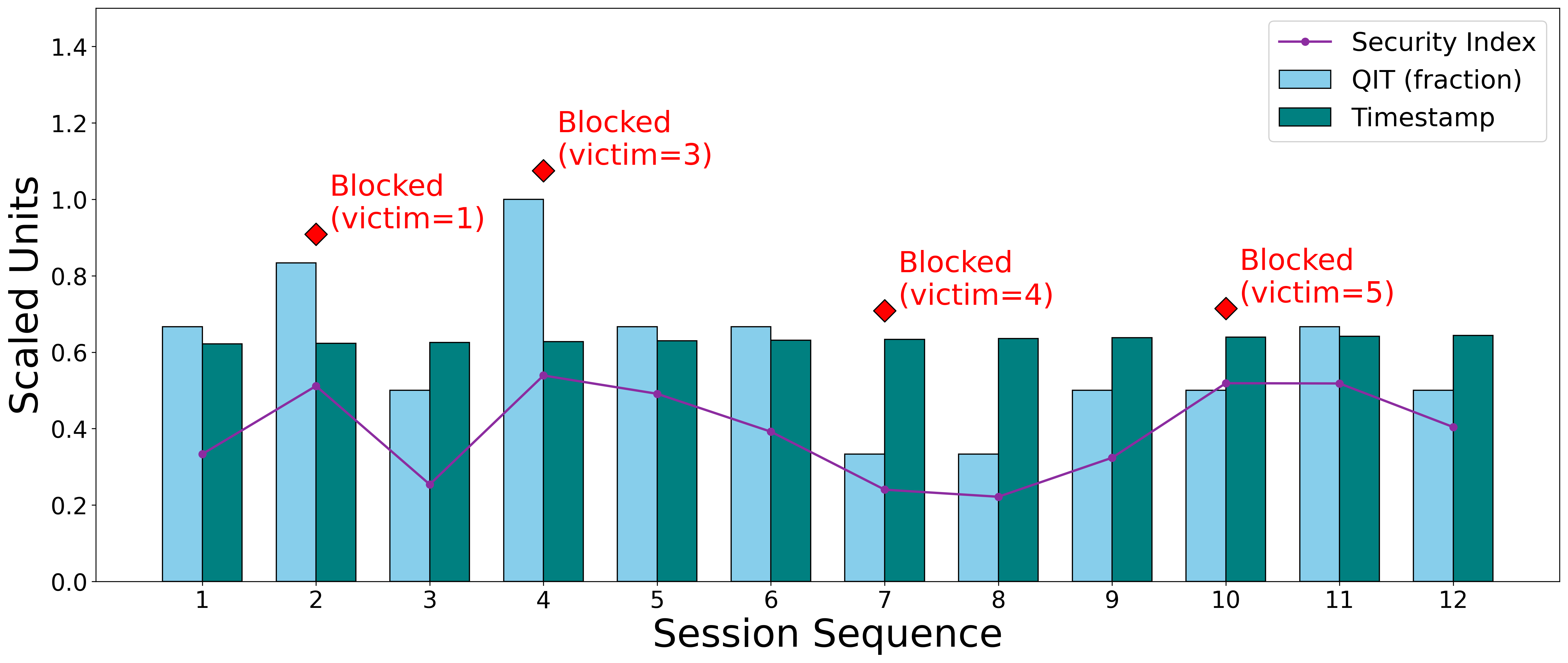}
         \caption{}
         \label{fig:}
     \end{subfigure}
        \caption{(a) and (b) Replay attacks with dynamic variations in QIT and timestamps.}
        \label{fig:replay attack}
\end{figure}
Fig.~\ref{fig:replay attack} presents a) the security performance of 25 sequential authentication sessions under a replay attack. The legitimate sessions between vehicles, nodes, and the gateway, plotted against simulated (logical) time. Most attempts were detected and rejected successfully, while only a few attacks succeeded. (b) detailed illustration of replay attack mitigation for individual sessions. Each session is characterized by a session-specific QIT and dynamic timestamp (T), collectively enforcing resistance to replay attacks (SG2) and forward and backward attacks (SG7). The combined security index, computed as a weighted sum of QIT strength and key variation between consecutive sessions, quantifies session robustness. Even sessions with relatively lower QIT or security index values remain protected because the uniqueness and freshness of session-specific parameters prevent accurate replay. Additionally, session key security (SG6) is maintained, as the dynamic combination of T and QIT values prevents adversaries from predicting previous or future session keys.

\begin{figure}[H]
\centering
\includegraphics[width=0.9\linewidth]{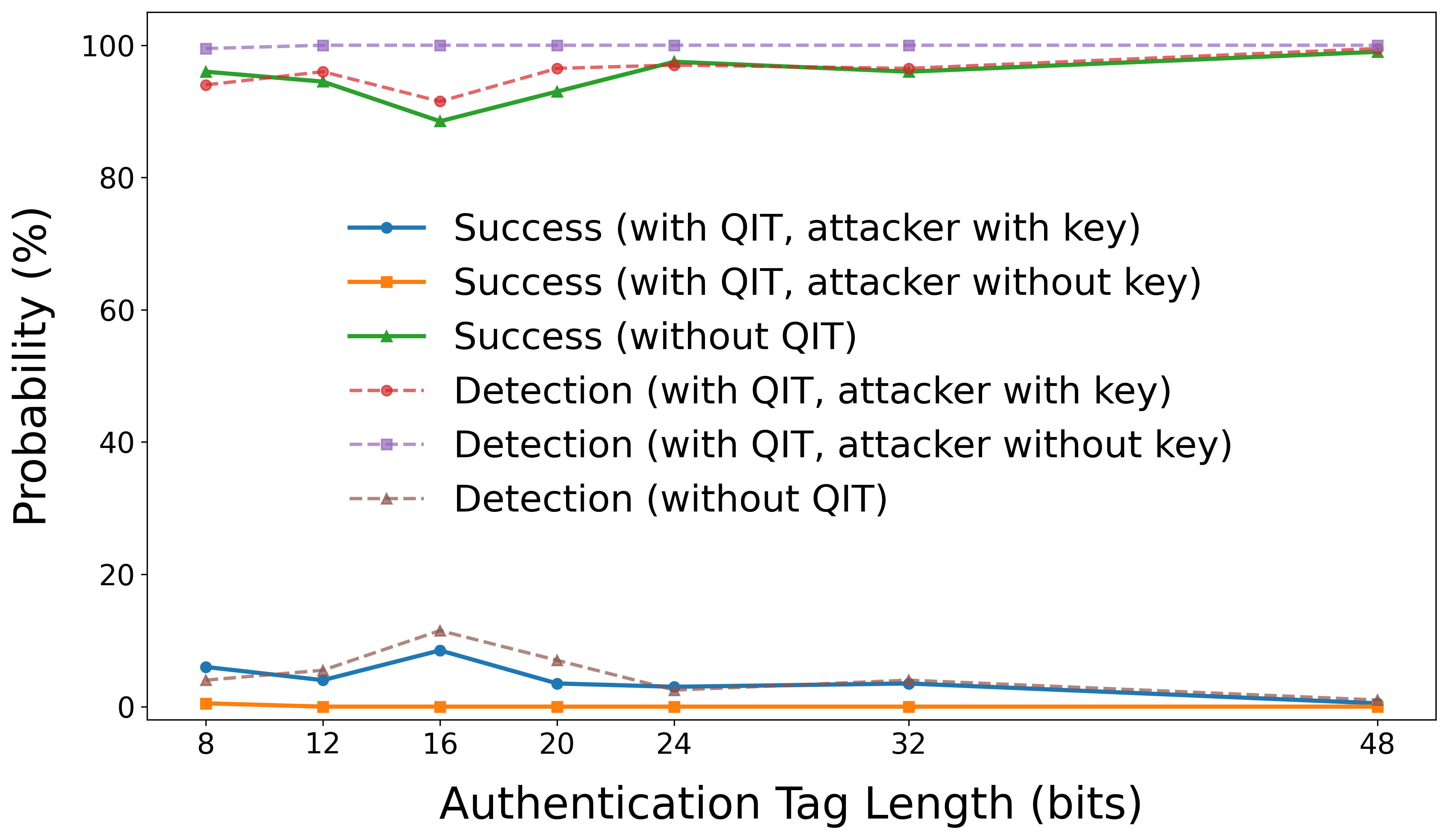}
\caption{MITM attack resistance under QIT authentication.}
\label{fig:mitm_attack}
\end{figure}

\begin{figure}[H]
\centering
\includegraphics[width=0.9\linewidth]{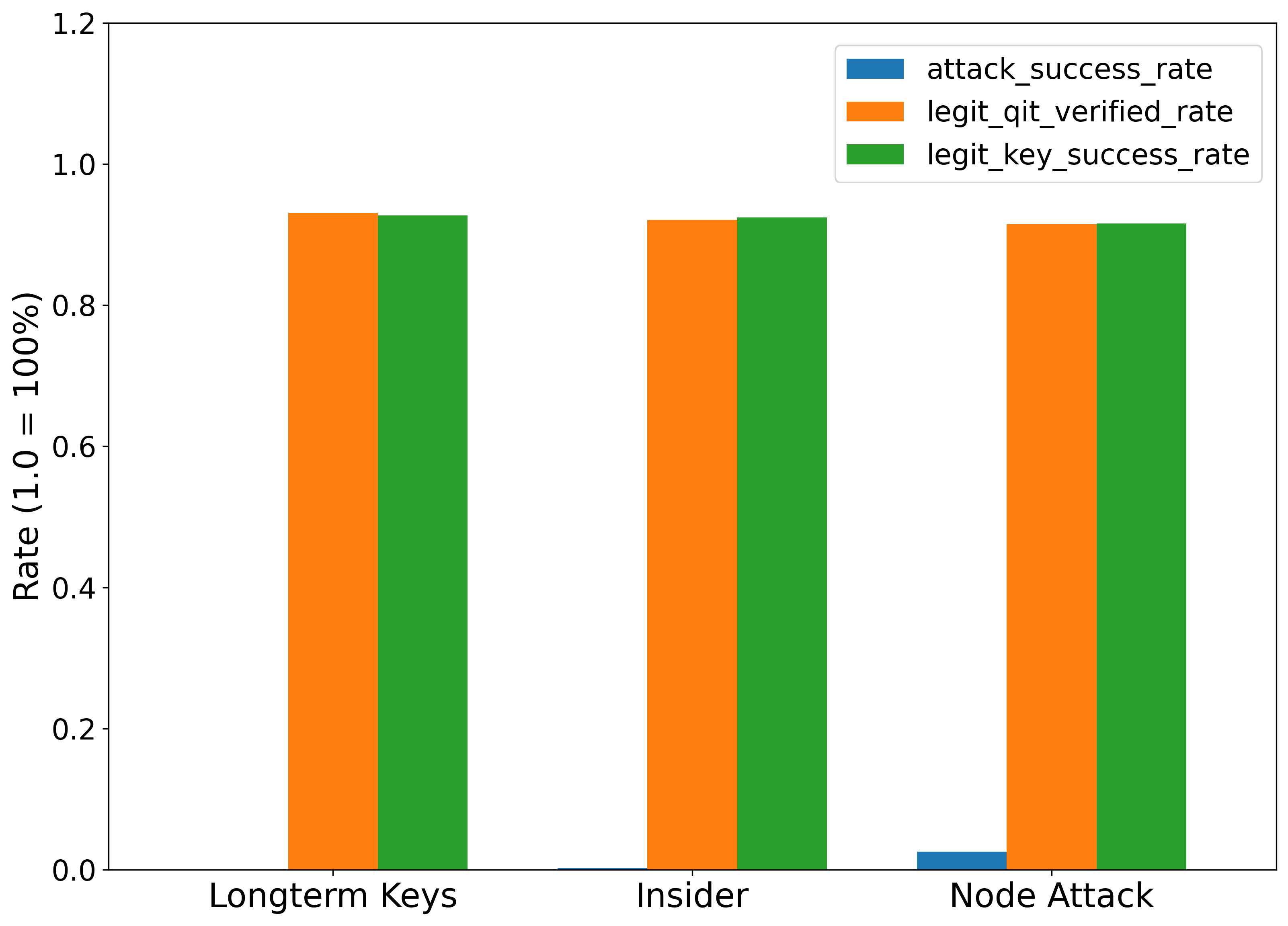}
\caption{Impersonation attacks}
\label{fig:impersonation_attack}
\end{figure}

Fig.~\ref{fig:mitm_attack} illustrates a simulation-based analysis demonstrating the resistance of the proposed scheme, QSAFE-V authentication framework, against MITM attacks. The success and detection probabilities are plotted as functions of the authentication tag length (in bits) for three different adversarial scenarios: (a) in a rare case, attackers possessing the legitimate key by guessing QIT successfully, (ii) in a highest condition, attackers unable to guess the the key under QIT-based authentication, and (iii) attackers in a purely classical setting without QIT. The results show that the probability of a successful forgery or message modification fractures due to authenticated quantum channels and session-specific parameters. Moreover, it significantly decreases as the tag length increases. This successfully represents the  Security Goal (SG3).
Fig.~\ref{fig:impersonation_attack} depicts simulated success rates for quantum identity authentication under three attack scenarios: long-term key compromise, insider attack, and node tampering or quantum eavesdropping. The attack success remains lower, while legitimate QIT verification and full key achieve the highest success rate. This demonstrates that the system provides high reliability for honest users while effectively resisting attacks, which highlights our proposed QSEFE-V model. 
\begin{figure}[H]
\centering
\includegraphics[width=0.80\linewidth]{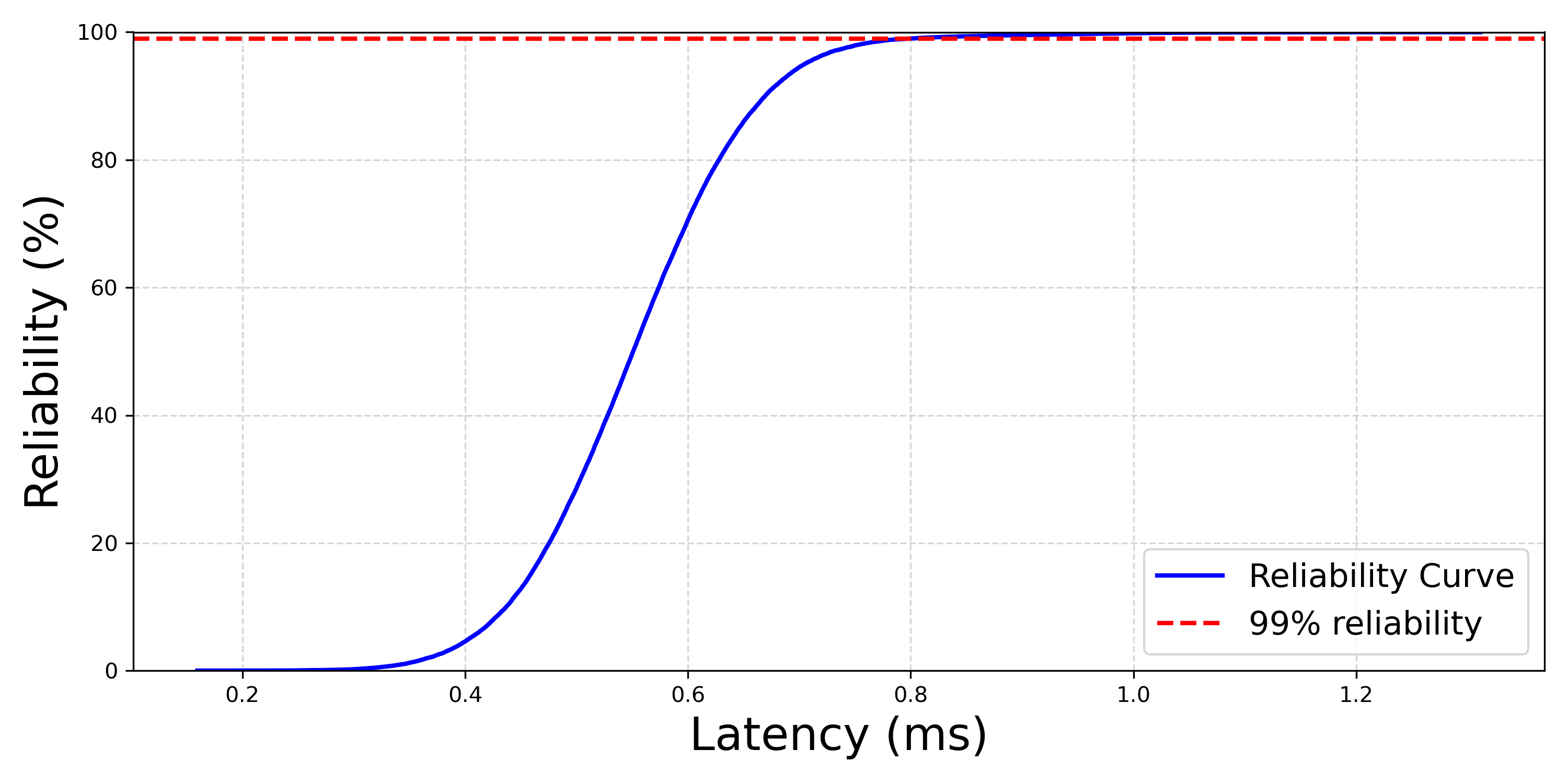}
\caption{Reliability vs Latency curve}
\label{fig:urllc}
\end{figure}
Fig.~\ref{fig:urllc} focuses on the simulation that evaluates end-to-end packet latency in the proposed QSAFE-V system, incorporating key processing stages, edge computation, QIT generation, verification, and configured grants, each with Gaussian jitter to reflect realistic variability. Local vehicular edge nodes reduce round-trip delays, enabling fast QIT verification and the highest reliability with the low latency value, which highlights the system’s capability under realistic Tactile Internet conditions.

\section{Conclusion}
This paper introduced QSAFE-V, a quantum-driven authentication framework designed specifically for edge-enabled vehicles in the Tactile Internet. By leveraging the principles of QKD and lightweight cryptographic mechanisms, QSAFE-V enables mutual authentication and secure key establishment with low computational overhead. The protocol is tailored to operate efficiently in latency-sensitive vehicular edge networks, addressing both classical and post-quantum attack surfaces.
Through detailed security analysis and comparative performance evaluations, we demonstrated that QSAFE-V achieves significant improvements in entropy strength, authentication latency, and attack resilience when compared to existing classical schemes such as IAR-AKA. Furthermore, QSAFE-V supports scalable deployment in vehicular edge environments while maintaining quantum-resilient properties.
Future work will involve extending QSAFE-V with adaptive quantum resource management, incorporating QML techniques for real-time threat detection, and validating the protocol in real-world vehicular edge testbeds using quantum simulators and hardware-based QKD modules.

\bibliographystyle{IEEEtran}
\bibliography{IEEEabrv,mybib}

\end{document}